\documentclass[11pt]{article}
\usepackage[T1]{fontenc}
\usepackage[utf8]{inputenc}
\usepackage{lmodern,microtype}
\usepackage{amsmath,amssymb,amsthm,booktabs,geometry}
\usepackage{xcolor}
\usepackage{url}
\usepackage{hyperref}
\hypersetup{colorlinks=true,linkcolor=blue!55!black,
 citecolor=blue!55!black,urlcolor=blue!55!black,
 pdftitle={A Computational Verification of the Markov Uniqueness Conjecture for z <= 10 to the power 500000},
 pdfauthor={Hu Tan and Ying Zhang},
 pdfsubject={Markov uniqueness, finite computation, and reproducibility}}
\newtheorem{theorem}{Theorem}[section]
\newtheorem{lemma}[theorem]{Lemma}
\newcommand{\N}{\mathbb{N}}
\newcommand{\B}{10^{500000}}
\newcommand{\code}[1]{\texttt{#1}}
\title{A Computational Verification of the Markov Uniqueness Conjecture for $z\leq 10^{500000}$}
\author{}
\date{}

\makeatletter
\renewcommand{\maketitle}{%
  \begin{center}

    \vspace*{-1.0em}

    {\LARGE
      \@title
      \par
    }

    \vspace{1.15em}

    {\large
      Hu Tan\textsuperscript{1}
      \hspace{2.0em}
      Ying Zhang\textsuperscript{2}
      \par
    }

    \vspace{0.85em}

    \begin{minipage}[t]{0.55\textwidth}
      \centering
      \small
      \mbox{\textsuperscript{1}Academy of Mathematics and Systems Science}\\
      Chinese Academy of Sciences\\
      Beijing 100190, China\\[0.25em]
      \footnotesize
      \href{mailto:tanhu2020@amss.ac.cn}
      {\texttt{tanhu2020@amss.ac.cn}}
    \end{minipage}%
    \hspace{0.035\textwidth}%
    \begin{minipage}[t]{0.40\textwidth}
      \centering
      \small
      \textsuperscript{2}School of Mathematical Sciences\\
      Soochow University\\
      Suzhou 215006, China\\[0.25em]
      \footnotesize
      \href{mailto:yzhang@suda.edu.cn}
      {\texttt{yzhang@suda.edu.cn}}
    \end{minipage}

    \vspace{0.9em}

  \end{center}

  \vspace{0.7em}
}
\makeatother

\begin{document}
\maketitle
\begin{abstract}
We describe a finite computational verification of Markov uniqueness
for normalized positive integer solutions of $x^2+y^2+z^2=3xyz$
with $z\leq10^{500000}$. The archived run reports
$239\,536\,376\,556$ triples with distinct maximal coordinates.
The method combines a disjoint tree decomposition, certified
fixed-point range bounds, a deterministic 48-bit candidate screen,
and replay modulo $2^{128}$. We give the mathematical justification
and a reproducibility procedure covering package checks,
retained-data auditing, and a full rerun. The conclusion depends on
correct execution of the exhaustive screen and replay; the retained
audit does not regenerate the deleted screening buckets.
\end{abstract}

\noindent\textit{2020 Mathematics Subject Classification:} Primary 11Y50; Secondary 11D25.

\noindent\textit{Keywords:} Markov equation; Markov number; Markov uniqueness conjecture;
computational verification; interval arithmetic.

\section{Introduction}
The positive integer solutions of the Markov equation
\begin{equation}\label{eq:markov}
 x^2+y^2+z^2=3xyz
\end{equation}
have a classical connection with Diophantine approximation and
indefinite binary quadratic forms. Both the spellings \emph{Markov}
and \emph{Markoff} occur in the literature. We use the former in the
text and preserve original spellings in references. A solution is
normalized when $1\leq x\leq y\leq z$; permutations are not counted
as different triples.

The uniqueness conjecture states that a positive integer is the largest
coordinate of at most one normalized Markov triple. It does not assert
that a Markov number occurs in only one triple in any coordinate.
Baragar \cite{baragar} gives this formulation, attributes the conjecture
to Frobenius, and reports a verification for maximal coordinates below
$10^{140}$. Button \cite{button} proves uniqueness when the largest
coordinate is prime. These results provide context; the present
computation does not depend on primality testing or the classification
of special families.

Zagier \cite{zagier} studies the counting function for normalized
triples and proves an asymptotic law with leading constant
$C=0.18071704711507\ldots$ multiplying $(\log X)^2$. His paper
also contains numerical data up to $10^{1300}$. Counting triples and
establishing uniqueness of their maxima are different tasks. We make
no priority claim and do not present these references as an exhaustive
survey of previous computational bounds.

This article documents the calculation archived on 26 October 2026
\cite{archive}. Its purpose is to state the finite result precisely,
expose the assumptions connecting the implementation to that result,
and provide reviewable evidence. New checks of the retained data are
reported separately from the original execution.

Write $\N=\{1,2,\ldots\}$. For $X\geq1$, define
\[
 \mathcal{T}(X)=\{(x,y,z)\in\N^3:x\leq y\leq z\leq X,
      \ x^2+y^2+z^2=3xyz\},\qquad N(X)=|\mathcal{T}(X)|.
\]
The tested statement is injectivity of
$\pi_X:\mathcal{T}(X)\to\N$, $(x,y,z)\mapsto z$.
The bound $\B$ is inclusive; 500000 is its exponent.

\medskip
\noindent\textbf{Reported computational finding.}
For $B=\B$, the archived implementation reports
\begin{equation}\label{eq:count}
 N(B)=239\,536\,376\,556
\end{equation}
and injectivity of $\pi_B$. The final record has matching visited
and expected counts and zero boundary ambiguities, missing candidate
groups, and 128-bit residue collisions. Subject to the correctness
conditions below, the number in \eqref{eq:count} is consequently also
the number of distinct maximal Markov numbers in the tested range.

Section~\ref{sec:method} gives the enumeration, range test, and collision
exclusion that make the computation possible. Section~\ref{sec:reproduce}
describes the retained evidence and the steps for checking it or
repeating the calculation.

\section{Method}\label{sec:method}

\subsection{Tree enumeration and partition}\label{sec:tree}
The Vieta involution replaces a coordinate $z$ by $3xy-z$.
After the exceptional triples $(1,1,1)$ and $(1,1,2)$, the normalized
tree starts at $(1,2,5)$ and has children
\begin{equation}\label{eq:children}
 L(a,b,c)=(a,c,3ac-b),\qquad
 R(a,b,c)=(b,c,3bc-a).
\end{equation}
For fixed $a\leq b$, the polynomial
$P(t)=t^2-3ab\,t+a^2+b^2$ has roots $c$ and $c'=3ab-c$.
Since $P(b)=a^2+(2-3a)b^2<0$ outside the exceptional cases,
$0<c'<b<c$. Mutating the maximum therefore gives the unique
descending neighbor; mutating either smaller coordinate gives an
ascending neighbor. Equal coordinates occur only in the exceptional
triples: if $a=b$, then $a\mid c$, and writing $c=ak$ gives
$k^2+2=3ak$, hence $k\mid2$ and $a=1$; the case $b=c$ is
excluded by $P(b)<0$. Descent decreases a positive integer and
terminates in the exceptional prefix. Reversing it proves that
\eqref{eq:children} enumerates every normalized triple exactly once.

Exact integer enumeration up to $H=2^{128}$ gives a prefix
$\mathcal P$ of 1456 nodes and 1455 frontier roots, each the first
node on its path with maximum exceeding $H$. Thus
\begin{equation}\label{eq:partition}
 \mathcal T(B)=\mathcal P\;\dot\cup\;
 \bigsqcup_{s=0}^{1454}\mathcal T_s(B),
\end{equation}
where each subtree includes its root. Increasing maxima and unique
descent give disjointness and coverage. Exact reconstruction of the
stored frontier checks the 1456 prefix records, the 1455 roots,
their ordering, and their 256-bit coordinate encodings. In particular,
every root has maximum below $2^{256}<B$.

For replay identification, $(1,2,5)$ is assigned $(u,v)=(1,1)$,
with $L:(u,v)\mapsto(u,u+v)$ and
$R:(u,v)\mapsto(u+v,v)$. Subtracting the smaller coordinate from
the larger recovers the path uniquely. The two exceptional prefix
nodes use $(0,1)$ and $(0,2)$. These labels identify nodes; they
are not the size bounds used next.

\subsection{Certified range decisions}\label{sec:bounds}
Large exact integers are avoided during the main traversal by storing
modular coordinates and bounds for
\begin{equation}\label{eq:f}
 f(t)=\operatorname{arcosh}(3t/2)
 =\log\frac{3t+\sqrt{9t^2-4}}2,\qquad F=f(B).
\end{equation}
Because $f$ is increasing, $z\leq B$ is equivalent to $f(z)\leq F$.

\begin{lemma}[Additive defect]\label{lem:defect}
For a normalized Markov triple $(a,b,c)$ with $c\geq5$,
\[
 0<f(a)+f(b)-f(c)<2/c^2.
\]
\end{lemma}
\begin{proof}
Put $A=3a$, $D=3b$, $C=3c$, and
$\Delta=A^2D^2-4(A^2+D^2)$. The larger-root relation
$2C=AD+\sqrt\Delta$ gives
\[
 \cosh(f(a)+f(b))-\cosh f(c)
 =\frac4{\sqrt{\Delta+16}+\sqrt\Delta}>0.
\]
Also $P(2b)=(5-6a)b^2+a^2<0$, so $c>2b$.
Together with $c'<b$, this gives
$\sqrt\Delta=3(c-c')>3c/2$; and
$\sinh f(c)=\sqrt{9c^2-4}/2>c$.
The mean value theorem bounds the defect by $4/(3c^2)<2/c^2$.
\end{proof}

Each frontier seed stores integer endpoints enclosing
$2^{255}f(x)/F$ and $2^{255}f(y)/F$.
Their containment is checked by an integer-only certificate,
independent of the high-precision calculation that proposed them.
At scale $2^{512}$, exact integer square roots give outward bounds
for the radicals in \eqref{eq:f}. After reducing a logarithm to
$1\leq m\leq2$, set $t=(m-1)/(m+1)\leq1/3$ and use
\[
 \log m=2\sum_{j=0}^{179}\frac{t^{2j+1}}{2j+1}+R,\qquad
 0\leq R\leq3(1/3)^{361}.
\]
Outward rounding and exact cross-multiplication then certify the
endpoints. For the denominator, use
$0<\log(3B)-f(B)<B^{-2}<2^{-512}$.
Indeed, with $s=4/(9B^2)$ and
$u=(1-\sqrt{1-s})/2$, this difference is $-\log(1-u)$,
bounded by $u/(1-u)\leq2/(7B^2)$.
All 2910 endpoint checks pass in the recorded certificate.

Let $q=2^{63}$. The certified Q255 endpoints $\ell_x,u_x$ are
converted to Q63 bounds by
\begin{equation}\label{eq:q63}
 l_x=\lfloor\ell_x/2^{192}\rfloor-1,\qquad
 h_x=\lceil u_x/2^{192}\rceil,
\end{equation}
and likewise for $y$. Bounds for each new maximum are obtained by
addition: $l_z=l_x+l_y$ and $h_z=h_x+h_y$.
This addition encloses the size function; the modular coordinates
continue to follow the Vieta recurrence.

To justify the lower bound, put $\delta_0=2/H^2=2^{-255}$.
Every maximum constructed in a frontier subtree exceeds $H$, so
its additive defect is below $\delta_0$.
The subtraction of one in \eqref{eq:q63} gives each smaller seed
coordinate a lower margin of at least $F/q>\delta_0$.
If both current coordinates have this margin, then
\[
 f(z)=f(x)+f(y)-\delta(z)
 \geq(F/q)(l_x+l_y)+\delta_0.
\]
The margin survives induction. Upper bounds follow from
$f(z)<f(x)+f(y)$. Consequently a child is accepted if $h_z\leq q$,
excluded if $l_z>q$, and otherwise marked ambiguous. Excluded
children may be pruned because all their descendants have larger
maxima. Zero ambiguities therefore gives an exhaustive range decision.

Roots are inserted before their children are tested, so their
preconditions are checked explicitly: positive lower bounds,
exact in-range root values, and root upper sums at most $q$.
For the supplied seeds, the smallest lower endpoint is
$7\,710\,275\,671\,045$ and the largest root upper sum is
$1\,423\,399\,665\,549\,849<q$.
At an included parent, positivity and $h_x+h_y\leq q$ imply
$2h_x+h_y<2q$ and $h_x+2h_y<2q$; hence child upper sums fit
in an unsigned 64-bit word. The implementation checks carries,
parent invariants, and path-label overflow. A positive overflow
in a lower sum certifies exclusion.

\subsection{Candidate screening and modular replay}\label{sec:screen}
The first traversal associates to a maximum $z$ the deterministic key
\[
 h(z)=\operatorname{Mix64}(z\bmod2^{64})\bmod2^{48}.
\]
The mixing function successively applies
$(t\mathbin{\mathrm{xor}}(t\gg30))\,\code{0xbf58476d1ce4e5b9}$,
$(t\mathbin{\mathrm{xor}}(t\gg27))\,\code{0x94d049bb133111eb}$,
and $t\mathbin{\mathrm{xor}}(t\gg31)$, with products modulo $2^{64}$.
The upper eight key bits select one of 256 buckets; the remaining
40 bits are stored in five bytes.

The subtrees are balanced into $G=13$ groups. Every unordered
group pair, including the diagonals, is processed, giving 91 jobs.
A job includes the prefix and either one group or two distinct
groups. Two nodes in the same group meet in its diagonal job;
two nodes in different groups meet in their joint job; the prefix
meets every group. Thus every repeated key appears as a duplicate
in at least one job. Each node is inserted only once within a job.

The bucket scanner stores the full 40-bit payload as
$\code{key}+1$ in an open-addressed table, with zero reserved for
an empty slot. Exact payload comparison and linear probing handle
table-address collisions; a load at most 0.60 leaves an empty slot.
The repeated keys from all jobs are sorted and deduplicated.
With $T=1456$ and $N=N(B)$, the total processed records are
\[
 G(N-T)+\frac{G(G+1)}2T
 =3\,113\,973\,008\,796,
\]
which includes repeated coverage across jobs.

A second traversal computes \eqref{eq:children} modulo $2^{128}$
and records every occurrence of a listed key. Each 28-byte record
contains a four-byte candidate index, a 16-byte residue, and two
four-byte path labels. The prefix and each frontier subtree have
separate shards. Polynomial evaluation commutes with reduction
modulo $2^k$, so exact seed residues determine the residues of
all descendants. The two-limb implementation propagates the
multiplication carry and subtraction borrow; equality compares
both 64-bit limbs.

The final merge requires at least two occurrences per candidate
key and checks for equal residues. It processes at most four
million keys per block, giving 26 blocks. A residue collision
requires exact investigation before any uniqueness conclusion.

\begin{theorem}[Collision exclusion]\label{thm:exclusion}
Suppose the traversal visits every element of $\mathcal T(B)$
exactly once, the screen lists every repeated key, and replay
records every occurrence with its correct residue modulo $2^{128}$.
If no two occurrences of any listed key have equal residues,
then the maximal-coordinate map $\pi_B$ is injective.
\end{theorem}
\begin{proof}
Equal maxima at distinct nodes give equal keys, so both occurrences
are replayed under a listed key. Their residues are equal,
contradicting the hypothesis.
\end{proof}
This implication is deterministic: unequal residues certify unequal
integers, regardless of the size of $B$. No statistical assumption
about the fingerprint distribution is required.

\subsection{Recorded result}\label{sec:results}
The final computation record in the supplementary materials
reports the quantities in Table~\ref{tab:main}. Its visited and
expected counts agree, and its three failure counters are zero.
The complete screen and replay, with the correctness conditions
above, yield the finite result stated in the introduction.

\begin{table}[htbp]
\centering
\caption{Recorded result for $z\leq10^{500000}$.}\label{tab:main}
\begin{tabular}{lr}
\toprule
Quantity & Value\\
\midrule
Visited nodes and expected nodes & $239\,536\,376\,556$\\
Prefix nodes / frontier subtrees & $1\,456\,/\,1\,455$\\
Groups / pair jobs & $13\,/\,91$\\
Candidate keys & $101\,871\,768$\\
Candidate occurrences & $203\,772\,186$\\
Replay shards / merge blocks & $1\,456\,/\,26$\\
Boundary ambiguities & $0$\\
Missing candidate groups & $0$\\
128-bit residue collisions & $0$\\
\bottomrule
\end{tabular}
\end{table}

\section{Reproducibility}\label{sec:reproduce}

\subsection{Supplementary materials and environment}
The supplementary materials accompanying this paper 
contain the computational evidence and reproducibility
instructions \cite{archive}. They were updated using Codex during
manuscript preparation. Consult these materials for the directory
layout, detailed commands, integrity checks, and test procedures.

The retained production evidence consists of the source programs,
frozen seeds and inputs, 91 candidate lists, merged candidates,
1456 replay shards, logs, manifests, and tests. The review scripts
and audit records used during manuscript preparation are described
below.

Production C\# code is compiled through PowerShell 7 with a .NET
runtime supporting \code{UInt128}; Python handles certificates,
metadata, grouping, and candidate merging. The supplied audit uses
Python and NumPy. The recorded review environment is Python
3.12.14, NumPy 2.3.5, and .NET 10.0.11; these are not established
as the original run's exact versions. The full production workflow
uses a Windows process-held lock. Use new working and scratch
directories, keep the archived evidence unchanged, and run the
original Python checks with assertions enabled.

\subsection{Verification and rerun procedure}
The following sequence separates checks on retained evidence from
a new execution. The supplementary instructions give the commands
and working-directory settings for each step.
\begin{enumerate}
\item \textbf{Check package integrity.}
Follow the integrity-checking procedure in the supplementary materials.
Check file lengths and hashes, frozen-input bindings, and the
target result fields. Hash agreement binds the files to the
manifest; it does not establish arithmetic correctness.

\item \textbf{Audit retained candidates and residues.}
Run the retained-data audit described in the supplementary materials
in a new working directory, allowing several GB of temporary space.
It reconstructs the exact frontier, checks all pair
memberships and the exact union of their saved candidate lists,
and validates every shard's structure, bindings, count, and digest.
It recomputes candidate fingerprints and multiplicities for all
$203\,772\,186$ occurrences. Independently of the production
merge, it partitions residues by their top six bits and sorts
each partition to detect equal 128-bit values.

The supplied audit record reports a pass on all 1456 shards,
zero repeated residues, and 4368 reconstructed sample paths
matching their stored residues. Candidate multiplicities are
$101\,843\,122$ keys occurring twice, $28\,642$ occurring three
times, and four occurring four times. Full residue reconstruction
is sampled; all saved residues are included in the equality test.

\item \textbf{Recompute certificates and small tests.}
Recompute the integer seed certificate and run the original exact
regression in a separate copy. At $10^{100}$, it compares 9670
exactly enumerated triples with an eight-bit screen that forces
collisions. Also run the supplementary tests for full-width
arithmetic, boundary decisions, and failure handling.
The supplied records report successful checks of all 16428
small-instance child decisions, a full-width test across 256
buckets, extremal modular arithmetic, classifier boundaries,
and malformed-input rejection. These tests exercise the code;
they do not replace the target computation.

\item \textbf{Recheck or regenerate replay.}
Follow the supplementary replay instructions to prepare a fresh
destination. Rechecking saved shards validates and merges the retained
records. A fresh traversal instead revisits all target nodes modulo
$2^{128}$; this step still uses the saved 48-bit candidate list.

\item \textbf{Repeat the full calculation.}
Follow the supplementary instructions for a full rerun in a separate
directory: regenerate and certify seeds, count the subtrees,
construct the 13-group plan, run all 91 pair jobs, merge candidates,
and execute fresh replay and its 26-block merge. The largest pair
has approximately 184.27 GB of five-byte payload; the supplied
guidance recommends at least 250 GB of available disk.
The original target replay used 16 workers. The archived
multi-day execution included interruptions and resumed jobs,
so its partial timers are not a reliable total-runtime estimate.
\end{enumerate}

A successful rerun at $B=10^{500000}$ must report the expected
$239\,536\,376\,556$ visited nodes, zero boundary ambiguities,
zero missing candidate groups, and zero 128-bit residue collisions,
with matching input and source bindings and a successful final status.
Agreement between counting and replay is a consistency check
because they share range logic.

The original temporary screening buckets were deleted and are
absent from the retained evidence. Their digests, pair candidate
lists, and replay shards remain. Accordingly, the retained-data
audit verifies the saved candidate union and occurrence records;
a fresh full calculation is required to regenerate and recheck
the original exhaustive screen. The review records document the
retained-data audit and small tests; a fresh target enumeration
was not performed during manuscript preparation. The result concerns the stated
finite range and does not resolve the unrestricted conjecture.

\subsection*{Declaration of AI use}
OpenAI's Codex was used extensively throughout this project,
including the development of the computational method; generation
and revision of the enumeration, screening, replay, and auditing
code; execution and inspection of computational runs; formulation
and checking of mathematical arguments; and drafting and editing
of this manuscript. Codex was also used to update and assemble the
supplementary reproducibility materials during manuscript preparation.
Both the original calculation and the
separately implemented retained-data audit were conducted within
this Codex-assisted workflow. The latter provides an additional
implementation check, rather than independent external replication.
The authors take responsibility for the mathematical claims,
software, computational evidence, and final text.

\end{document}